\RequirePackage{fix-cm}
\documentclass[envcountsame,envcountsect,smallextended]{svjour3}       
\smartqed  
\usepackage{graphicx}
\usepackage[latin1]{inputenc}
\usepackage{longtable,helvet,courier,fancyhdr}
\usepackage{amsfonts}
\usepackage{amsmath,amssymb,bbm}
\usepackage{multicol}
\usepackage{algorithm}
\usepackage{algorithmic}

\def\cocoa{{\hbox{\rm C\kern-.13em o\kern-.07em C\kern-.13em o\kern-.15em A}}}

\def\H{\mathcal{H}}
\def\I{\mathcal{I}}

\def\J{\mathcal{J}}

\def\mf{\mathfrak{m}}
\def\mult#1#2#3{{#1}_{#2}(#3)}

\def\P{\mathcal{P}}
\def\QQ{\mathbbm{Q}}

\def\X{\mathcal{X}}
\def\B{\mathcal{B}}
\def\I{\mathcal{I}}
\def\K{\mathcal{K}}
\def\J{\mathcal{J}}

\def\kk{\mathbbm{k}}
\def\ZZ{\mathbbm{Z}}
\def\mf{\mathfrak{m}}

\def\P{\mathcal{P}}

\def\R{\mathcal{R}}
\def \LT{{\rm LT}}
\def \LM{{\rm LM}}
\def \LC{{\rm LC}}
\def \HF{{\rm HF}}

\def \HP{{\rm HP}}
\def \HS{{\rm HS}}
\def \N{{\rm N}}
\def \hilb{{\rm hilb}}

\def \VF{{\rm VF}}
\def \VP{{\rm VP}}
\def\ri{\rangle}
\def\li{\langle}

\newcommand{\leftexp}[2]{{\vphantom{#2}}^{#1}{#2}}

\DeclareMathOperator{\cls}{cls}

\DeclareMathOperator{\depth}{depth}

\DeclareMathOperator{\reg}{reg}
\DeclareMathOperator{\sat}{sat}

\journalname{AAECC}
\begin{document}

\title{A Pommaret Bases Approach to the Degree of a Polynomial Ideal}
\author{Bentolhoda Binaei \and Amir Hashemi \and Werner M. Seiler
}

\authorrunning{Binaei et al.} 

\institute{Bentolhoda Binaei \at
              Department of Mathematical Sciences, 
              Isfahan University of Technology\\ Isfahan, 84156-83111, Iran\\
              \email{H.Binaei@math.iut.ac.ir}
           \and
           Amir Hashemi \at
           Department of Mathematical Sciences, 
           Isfahan University of Technology\\ Isfahan, 84156-83111, Iran;\\
\email{Amir.Hashemi@cc.iut.ac.ir}
\and Werner M. Seiler \at
Institut f\"{u}r Mathematik,
Universit\"at Kassel\\
Heinrich-Plett-Stra\ss e 40, 34132 Kassel, Germany\\
\email{Seiler@mathematik.uni-kassel.de}
}

\date{Received: / Accepted: date}

\maketitle

\begin{abstract}
  In this paper, we study first the relationship between Pommaret bases and Hilbert
  series.  Given a finite Pommaret basis, we derive new explicit formulas
  for the Hilbert series and for the degree of the ideal generated by it
  which exhibit more clearly the influence of each generator.  Then we
  establish a new dimension depending B\'ezout bound for the degree and use
  it to obtain a dimension depending bound for the ideal membership
  problem.

\end{abstract}

\section{Introduction}\label{sec:1}

Gr\"obner bases, introduced by Bruno Buchberger in his PhD thesis (see
\cite{Bruno1,Bruno2}), have become a powerful tool for constructive
problems in polynomial ideal theory and related domains.  It is well-known
that they allow us to design algorithms for computing important invariants
like dimension, degree and Hilbert function. However, the bases themselves
are largely independent of the values of these invariants. This is in
marked contrast to Pommaret bases which reflect many combinatorial and
homological properties of the ideals they generate. They are not only of
computational interest by providing easy access to many invariants relevant
for algebraic geometry, but also allow for alternative constructive proofs
of many theoretical results and thus lead to a much closer intertwining of
computation and theory than ordinary Gr\"obner bases.

Pommaret bases are a particular form of involutive bases which in turn
represent a special kind of Gr\"obner bases with additional combinatorial
properties.  Involutive bases combine ideas of the Janet-Riquier theory of
partial differential equations with Gr\"obner bases.  Zharkov and Blinkov
\cite{zharkov} introduced Pommaret bases as {\em involutive polynomial
  bases} into commutative algebra. Later, Gerdt, Blinkov, Zharkov, and
others developed a general theory of involutive bases \cite{Blinkov}.  The
terminology Pommaret bases is historically incorrect, as they appear
already in the work of Janet (see e.\,g.\ \cite[pp. 30--31]{Janet});
however, the name has become standard by now.  For a general survey on
involutive bases with special emphasis on Pommaret bases see e.\,g.\
\cite{wms:comb1,wms:comb2,wms:invol} (the last reference also contains some
historical notes).

In the sequel, we will apply some of the above mentioned results to
effective algebraic geometry.  {\em B\'ezout's theorem} may be considered
as a generalization of the fundamental theorem of algebra.  Let
$f_1,\ldots ,f_{n-1}\in \P=\kk[x_1,\ldots ,x_n]$ be a sequence of
homogeneous polynomials. If the system $f_1=\cdots=f_{n-1}=0$ has a finite
number of projective solutions (i.\,e.\ the dimension as projective variety
is zero), then the number of solutions, counted with multiplicities, is at
most $\prod_{i=1}^{n-1}\deg(f_i)$, see \cite[page 174]{Kemper}. One can
consider higher-dimensional extensions of this result using the {\em
  degree} of an ideal. If the ideal $\I$ is generated by $k$ homogeneous
polynomials of degrees $d_1\ge \cdots \ge d_k$, then
$\deg{(\I)}\le d_1\cdots d_\mu$ with $\mu=\min\{k,n\}$, see e.\,g.\
\cite[Lem.~3]{Brown2}, \cite[Thm.~1]{Heintz} or \cite[Lem.~2.95]{Ritscher}.
We will refer to this upper bound as {\em B\'ezout's bound}.

In this article, we are mainly concerned with two related topics.  Firstly,
we will express the degree of an ideal in terms of the degrees and classes
of the elements of its Pommaret basis.  As a by-product, we will provide a
new proof for the rationality of the Hilbert series and an explicit formula
for its numerator. Secondly, we will derive the new {\em dimension
  depending} B\'ezout bound $d_1\cdots d_{n-D}$ where $D=\dim(\I)$ for the
degree of $\I$.  Masser and W\"ustholz \cite[Thm.~II]{Masser} proposed the
upper bound $d_1^{n-D}$ (see also \cite[Prop.~3.5]{bayer3} and
\cite{Brown2}). However, our bound is sharper and our proof is more
elementary. While the bound itself is independent of Pommaret bases, our
proof relies crucially on special properties available only in quasi stable
position, i.\,e.\ the generic position characterized by the existence of a
finite Pommaret basis.

The article is organized as follows. In the next section, we review the
basic definitions and notations which will be used throughout. Section 3
investigates the relationship between Pommaret bases and Hilbert series. In
the last section, we derive our dimension depending B\'ezout bound and
bound for the ideal membership problem, respectively.

\section{Preliminaries}
\label{sec:2}

We first introduce basic notations and preliminaries needed in the
subsequent sections. Throughout $\P=\kk[x_{1},\dots,x_{n}]$ will be the
polynomial ring over an infinite field $\kk$. We consider always
homogeneous polynomials $f_1,\ldots ,f_k\in \P$ and the ideal
$\I=\li f_1,\ldots ,f_k \ri$ generated by them.  We assume that each $f_i$
is non-zero and denote its total degree by $d_i$.  We sort the $f_i$ so
that $d_1\ge d_2\ge \cdots\ge d_k$ and we set $d=d_1$.  Furthermore, the
dimension of $\I$, denoted by $D=\dim(\I)$, is the Krull
dimension\footnote{Note that the Krull dimension corresponds to the
  dimension as \emph{affine} and not as \emph{projective} variety, although
  we work exclusively with homogeneous ideals.  We stick with the affine
  picture to facilitate comparison with other results which are also based
  on the dimension as affine variety.} of the corresponding factor ring
$\R=\P/\I$.  It is trivial that $k\ge n-D$. If $f\in \P$, the equivalence
class of $f$ w.r.t. the congruence relation modulo $\I$ is an element of
$\R$ and is denoted by $[f]$.  Finally, we work throughout with the degree
reverse lexicographic term order with $x_n\prec \cdots \prec x_1$.

The {\em leading term} of a polynomial $0\neq f\in \P$, denoted by
${\LT}(f)$, is the greatest term appearing in $f$ and its coefficient is
the {\em leading coefficient}, denoted by ${\LC}(f)$; the {\em leading
  monomial} is the product ${\LM}(f)={\LC}(f){\LT}(f)$. The {\em leading
  ideal} of $\I$ is the monomial ideal
${\LT}(\I) = \langle {\LT}(f) \ | \ 0\neq f \in {\I}\rangle$. A finite
subset $\{g_1, \ldots , g_m\} \subset {\I}$ is called a {\em Gr\"obner
  basis} of $\I$ for $\prec$, if
${\LT}(\I) = \langle \LT(G) \rangle=\langle \LT(g)\mid g\in G\rangle$.  We
refer to \cite{Becker,little} for more details on Gr\"obner bases.  We will
denote by $\deg (\I, \prec)$ the maximal degree of the elements of the
reduced Gr\"obner basis of $\I$ with respect to $\prec$ (see
\cite{Lazard81,Daniel83,Giusti1}).

For a positive integer $s$, we denote by $\R_{s}$ the set of elements of
the factor ring $\R$ of degree $s$. Then the {\em Hilbert function} of $\I$
is defined by $\HF_{\I}(t)=\dim_{\kk}(\R_{t})$ where $\dim_{\kk}$ denotes
the dimension as a $\kk$-vector space. From a certain degree on, this
function of $t$ is equal to a (unique) polynomial in $t$, called {\em
  Hilbert polynomial}, and denoted by $\HP_{\I}$. The {\em Hilbert
  regularity} of $\I$ is
$$\hilb(\I)=\min \{m  \ \arrowvert \  \forall t \geq m , \
\HF_{\I}(t)=\HP_{\I}(t)\}\,.$$ We have the identity
$\dim(\I)=\deg(\HP_{\I})+1$, see \cite[Thm.~12, page 464]{little} and by
Macaulay's theorem $\HF_{\I}=\HF_{\LT(\I)}$.

The {\em Hilbert series} of $\I$ is the power series
${\HS_{\I}}(t)=\sum_{s=0}^{\infty}{{\HF_{\I}}(s)t^{s}}$.  This series can
be expressed as a quotient $\HS_\I(t)=N(t)/(1-t)^D$ with a polynomial
$N\in\QQ[t]$ satisfying $N(1)\ne 0$ (see \cite[Thm.~7, page 130]{Ralf} or
\cite{Valla}).  In the next section, we will provide a new proof of this
fact using Pommaret bases.

\begin{definition}[{\cite[page 52]{Hart}}]
  If $D>0$, then the degree of $\I$, denoted by $\deg(\I)$,\footnote{Please
    note that despite the similarity in notation $\deg{(\I)}$ and
    $\deg{(\I,\prec)}$ refer to very different objects!} is $(D-1)!$ times
  the leading coefficient of the Hilbert polynomial of $\I$. If $D=0$, then
  $\deg(\I)$ is defined to be the sum of the coefficients of $\HS_{\I}(t)$.
\end{definition}

\begin{remark}\label{rem:lendeg}
  By \cite[page 173]{Kemper}, we have $\deg(\I)=N(1)$ and in consequence since $\I$ and $\LT(\I)$ share the same Hilbert function, $\deg(\I)=\deg(\LT(\I))$.  We also need the
  following auxiliary results on the degree of an ideal from \cite{Brown2}.
  Let $Q$ be a $\mathfrak{p}$-primary ideal. We say that $Q$ has {\em
    length} $\ell$, if there exists a chain
  \begin{displaymath}
    Q=Q_1\subset \cdots \subset Q_\ell=\mathfrak{p}
  \end{displaymath}
  of primary ideals $Q_1,\ldots ,Q_\ell$, but no longer chain of this form.
  Let $r$ be the least positive integer $a$ such that
  $\mathfrak{p}^a\subset Q$.  With these
  notations, we find that $r\le \ell$ and
  $\deg(Q)=\ell\deg(\mathfrak{p})$, cf.~\cite[page 282]{Brown2}. Furthermore, the degree of $\I$ is
  equal to the sum of the degrees of its primary components of dimension
  $D$.
\end{remark}

Let $\mathfrak{m}=\li x_1,\ldots ,x_n\ri$ be the unique homogeneous maximal
ideal of $\P$.  The ideal $\I^{\sat}=\I:\mathfrak{m}^\infty$ is called the
{\em saturation} of $\I$.  The {\em satiety} of $\I$, denoted by
$\sat(\I)$, is the smallest positive integer $m$ such that
$\I_{\ell}=\I^{\sat }_{\ell}$ for all $\ell\ge m$.  It is always a finite
number \cite[Rem.~1.3]{bayer_stillman}.  As a trivial consequence, $\I$ and
$\I^{\sat}$ possess the same Hilbert polynomial and thus in particular the
same degree.  By \cite[Lem.~1.6]{bayer_stillman}, $\I^{\sat}=\I:y^{\infty}$
for a generic linear form $y$.

\begin{definition}\label{defbayer}
  The ideal $\I$ is {\em $m$-regular}, if its $i$th syzygy module can be
  generated by elements of degree at most $m+i$.  The {\em
    Castelnuovo-Mumford regularity} $\reg(\I)$ is the smallest $m$ such
  that $\I$ is $m$-regular.
\end{definition}

For more details on the regularity, we refer to
\cite{bayer_stillman,bermejo}.  The polynomials $f_1,\ldots ,f_k\in \P$
form an $\I$-\emph{regular sequence} for an ideal $\I$, if they generate a
proper ideal in $\P$ and if $[f_j]$ is a non zero divisor in the ring
$\P/\li \I,f_1,\ldots ,f_{j-1}\ri$ for $j=1,\ldots ,k$.  We simply speak of
a \emph{regular sequence}, if $\I=0$.  The {\em depth} of $\I$, denoted by
$\depth(\I)$ is the maximal length of an $\I$-regular sequence.

Given a polynomial $f\in \P$ with $\LT(f)=x^\alpha$ where
$\alpha=(\alpha_1,\ldots ,\alpha_n)$, the \emph{class} of $f$ is the
integer $\cls({f}) = \max{\{i\mid\alpha_{i}\neq0\}}$. Then the
\emph{multiplicative variables} of $f$ are
$\mult{\X}{P}{f}=\{x_{\cls({f})},\ldots ,x_n\}$.\footnote{We follow here
  the conventions of \cite{Blinkov}.  In \cite{wms:comb2}, a convention is
  used which corresponds to reverting the order of the variables
  $x_{1},\dots,x_{n}$.  This implies e.\,g.\ that the class is defined as
  the minimum and not the maximum.  Thus care must be taken when
  transferring results of \cite{wms:comb2} to the conventions used in this
  article.} The term $x^{\beta}$ is a \emph{Pommaret divisor} of
$x^{\alpha}$, written $x^{\beta}\mid_P x^{\alpha}$, if
$x^{\beta}\mid x^{\alpha}$ and
$x^{\alpha-\beta}\in\kk[\mult{\X}{P}{x^{\beta}}]$.

\begin{definition}\label{def:pombas}
  Let $\H\subset \I$ be a finite set such that no leading term of an
  element of $\H$ is a Pommaret divisor of the leading term of another
  element. Then $\H$ is a \emph{Pommaret basis} of $\I$, if
  \begin{displaymath}
    \bigoplus_{h\in\H}\kk[\mult{\X}{P}{h}]\cdot\LT(h)=\LT(\I)\,.
  \end{displaymath}
\end{definition}

If an ideal $\I$ has a Pommaret basis $\H$, then $\reg(\I)$ equals the
maximal degree of an element of $\H$ and $\depth(\I)$ is given by $n$ minus
the maximal class of an element of $\H$.  Furthermore
$\I^{\sat}=\I:x_{n}^{\infty}$ and the satiety is the maximal degree of an
element of class $n$ in $\H$.  We refer the reader to
\cite{wms:comb1,wms:comb2} and \cite[Chap.~3-5]{wms:invol} for a thorough
introduction into the theory of Pommaret bases.

It follows immediately from the definition that any Pommaret basis is a
(generally non-reduced) Gr\"obner basis.  The main difference between
Gr\"obner and Pommaret bases lies in the fact that any polynomial $f\in \I$
has a {\em unique} involutive standard representation, i.\,e.\ a standard
representation where all coefficients depend only on the multiplicative
variables of the respective generator, by the following direct sum
decomposition as graded $\kk$-linear spaces:
\begin{equation}\label{eq:pombas}
  \bigoplus_{h\in\H}\kk[\mult{\X}{P}{h}]\cdot h=\I\,.
\end{equation}
It furthermore allows to read off immediately the \emph{volume function} of
$\I$
\begin{equation}\label{eq:vf}
  \VF_{\I}(t)=\dim_{\kk}(\I_{t})=\sum_{h\in\H}[\deg{(h)}\leq
  t]\binom{t-\deg{(h)}+|\mult{\X}{P}{h}|-1}{t-\deg(h)}
\end{equation}
where $[\cdot]$ denotes the Kronecker-Iverson symbol which yields $1$, if
the condition in the bracket is satisfied, and $0$ otherwise.  Obviously,
the volume function is closely related to the Hilbert function:
$\HF_{\I}=\VF_{\P}-\VF_{\I}$.  Thus we obtain without any further computation
the Hilbert function from a Pommaret basis.  The same is true of the
Hilbert polynomial: the volume polynomial $\VP_{\I}$ is given by the same
expression as $\VF_{\I}$ with only the Kronecker-Iverson symbol omitted and
then $\HP_{\I}=\VP_{\P}-\VP_{\I}$.

While this observation allows for an easy computation of both Hilbert
function and polynomial for any concrete ideal possessing a Pommaret basis,
it is not satisfying from a theoretical point of view.  Via
\eqref{eq:pombas} we have only easy access to the volume function; the
dependence of properties of the Hilbert function and related invariants
like the degree on properties of the Pommaret basis is more difficult to
assess.  Therefore, we will exhibit this relationship in more detail in the
next section.\footnote{In \cite{wms:comb2,wms:invol} also
  \emph{complementary decompositions}, i.\,e.\ direct sum decompositions of
  the complement of $\LT(\I)$ are discussed and it is shown that any
  Pommaret basis induces one.  Then one can write down an explicit formula
  for $\HF_{\I}$ with a similar structure as \eqref{eq:vf}.  However, this
  only transforms the problem into understanding the precise relationship
  between the complementary decomposition and the Pommaret basis.  While
  this is relatively simple with regard to, say, $\dim(\I)$ and
  $\depth(\I)$ (see the corresponding results in
  \cite{wms:comb2,wms:invol}), the situation becomes non-trivial for
  $\deg(\I)$.}

Unfortunately, Pommaret bases do not always exist. However, this is only a
question of the used variables: since we assume that $\kk$ is an infinite
field, any ideal possesses a Pommaret basis after a generic linear change
of variables \cite{wms:comb2}.  More precisely, we meet here the
combinatorial notion of quasi-stability.\footnote{Quasi stable ideals are
  also known by many other names like weakly stable ideals, ideals of
  nested type or ideals of Borel type.}

\begin{definition}
  A monomial ideal $\J$ is called \emph{quasi stable}, if for any term
  $m \in \J$ and all integers $i, j,s$ with $1 \le j < i \le n$ and $s>0$,
  if $x_i^s\mid m$ there exists an integer $t\ge 0$ such that
  $x_j^tm/x_i^s\in \J$. The polynomial ideal $\I$ is in {\em quasi stable
    position} if $\LT(\I)$ is quasi stable.
\end{definition}

\begin{proposition}[{\cite{wms:comb2}}]
  A monomial ideal $\J$ has a Pommaret basis, if and only if it is quasi
  stable.  A polynomial ideal $\I$ has a Pommaret basis, if and only if it
  is in quasi stable position.
\end{proposition}

\begin{remark}\label{linchen}
  It is trivial to see that all the objects studied in this work like the
  Hilbert function or the regularity remain invariant under linear changes
  of coordinates.  Hence we may in the sequel always assume that we are in
  quasi stable position and thus that $\I$ has a Pommaret basis $\H$.
\end{remark}

In the sequel, we will use the following notations: given an ideal $\I$ in
quasi stable position, we write $\H=\{h_1,\ldots ,h_s\}$ for its Pommaret
basis. Furthermore, for each $i$ we set $m_i=\LT(h_i)$, $c_i=\deg(m_i)$ and
$n_i=n-|\mult{\X}{P}{m_i}|$, the number of non-multiplicative variables of
$m_i$. By definition, $\mult{\X}{P}{h_i}=\mult{\X}{P}{m_i}$ and hence
$\{m_1,\ldots ,m_s\}$ forms a Pommaret basis for $\LT(\I)$.

\section{Pommaret Bases and Hilbert Series}
\label{sec:8}

We now study the relationship between the Pommaret basis $\H$ of a
polynomial ideal $\I$ and certain invariants of $\I$ related to its Hilbert
series. Our main results are new formulae expressing $\deg{(\I)}$ and the
coefficients of the numerator $N$ of the Hilbert series in terms of the
degrees $c_{i}$ and the numbers $n_{i}$ of non-multiplicative variables of
the elements $h_{i}$ of $\H$. As a by-product, we provide new proofs of
some classical results like the following one.

\begin{theorem}
  \label{th1sec8}
  The Hilbert series can be written as a rational function of the form
  $\HS_\I(t)=N(t)/(1-t)^D$ with a polynomial $N\in\ZZ[t]$ satisfying
  $N(1)\neq0$.
\end{theorem}

\begin{proof}
  By Rem. \ref{linchen}, we may assume that $\I$ is in quasi stable
  position.  It is easy to see that the Taylor coefficient of order $q$ of
  the function $1/(1-t)^{n}$ and the number of terms of degree $q$ in $n$
  variables coincide and thus the Hilbert series of the full polynomial
  ring is given by $1/(1-t)^{n}$.  It then follows from \eqref{eq:vf} and
  the subsequent discussion that
  \begin{equation}\label{eq:hs}
    \HS_{\I}(t)=\frac{1}{(1-t)^n}-
                \sum_{i=1}^{s}\frac{t^{\deg(m_i)}}{(1-t)^{|\mult{\X}{P}{m_i}|}}
               =\frac{1-\sum_{i=1}^{s}{(1-t)^{n_i}t^{c_i}}}{(1-t)^n}\,.
  \end{equation}
  We number the generators $h_{i}$ such that
  $m_1,\ldots, m_\ell\in \kk[x_1,\ldots ,x_{n-D}]$ and
  $m_{\ell+1},\ldots ,m_s\notin \kk[x_1,\ldots ,x_{n-D}]$. Since $\I$ is in
  quasi stable position, the set $\B=\{m_1,\ldots ,m_\ell\}$ contains pure
  powers of the variables $x_1,\ldots ,x_{n-D}$.  Hence, if we consider the
  contraction $\J=\LT(\I)\cap\P'$ with $\P'=\kk[x_1,\ldots ,x_{n-D}]$, then
  $\B$ is its Pommaret basis which trivially implies that $\J$ is
  zero-dimensional. It follows that the Hilbert series of $\J$ is a
  polynomial $P(t)$ with $P(1)\ne 0$ and the Hilbert series of the
  extension ideal $\J^{e}=\langle\J\rangle_{\P}\subset\P$ is given by
  $P(t)/(1-t)^D$.

  On the other hand, the Pommaret basis $\B$ induces the direct sum
  decomposition
  $\J^{e}=\bigoplus_{i=1}^{\ell}\kk[\mult{\X}{P}{m_i}]\cdot m_i$ and hence
  \begin{displaymath}
    \frac{P(t)}{(1-t)^D}=\frac{1-\sum_{i=1}^{\ell}{(1-t)^{n_i}t^{c_i}}}
                              {(1-t)^n}\,.
  \end{displaymath}
  Consequently,
  $1-\sum_{i=1}^{\ell}{(1-t)^{n_i}t^{c_i}}=(1-t)^{n-D}P(t)$. Our
  assumptions imply that for $i=\ell+1,\ldots,s$ the term $m_i$ contains at
  least one of the variables $x_{n-D+1},\ldots ,x_n$ entailing that all the
  variables $x_1,\ldots,x_{n-D}$ are non-multiplicative for it and thus
  $n_i\ge n-D$. It follows that the Hilbert series of $\I$ can be expressed
  in the form
  \begin{displaymath}
    \frac{(1-t)^{n-D}P(t)-\sum_{i=\ell+1}^{s}{(1-t)^{n_i}t^{c_i}}}{(1-t)^n}=
    \frac{P(t)-\sum_{i=\ell+1}^{s}{(1-t)^{n_i-n+D}t^{c_i}}}{(1-t)^{D}}\,.
  \end{displaymath}
  Writing $N(t)$ for the numerator of the quotient on the right hand side,
  we thus obtain the rational form $\HS_\I(t)=N(t)/(1-t)^D$.

  There only remains to show that $N(1)\ne 0$. For $n_i-n+D>0$ the
  polynomial $(1-t)^{n_i-n+D}t^{c_i}$ vanishes at $t=1$. Thus it suffices
  to consider only $P(t)$ (a polynomial with non-negative integer
  coefficients) and the polynomials $-(1-t)^{n_i-n+D}t^{c_i}$ for those
  indices $\ell+1\leq i\leq s$ with $n_i-n+D=0$ (which at $t=1$ yields
  $-1$). The condition $n_i-n+D=0$ corresponds to a leading term of the
  form
  $m_i=x_1^{\alpha_1}\cdots
  x_{n-D}^{\alpha_{n-D}}x_{n-D+1}^{\alpha_{n-D+1}}$ with
  $\alpha_{n-D+1}\ne 0$. Since $\{m_1,\ldots ,m_s\}$ forms a Pommaret
  basis, $x_1^{\alpha_1}\cdots x_{n-D}^{\alpha_{n-D}}\notin \J$ and
  therefore we can associate the unique non-zero class
  $[x_1^{\alpha_1}\cdots x_{n-D}^{\alpha_{n-D}}]$ in the quotient ring
  $\P'/\J$ to each $m_i$ with $n_i-n+D=0$. Since $\dim(\I)=D$, there is no
  pure power of $x_{n-D+1}$ contained in $\I$ and thus we cannot find any
  $m_i$ with $n_i-n+D=0$ which would correspond to the class
  $[1]$. Therefore, the number of leading terms $m_i$ with $n_i-n+D=0$ is
  strictly less than $\dim_\kk(\P'/\J)=P(1)$ and $N(1)\ne 0$. \qed
\end{proof}

\begin{remark}
  The basic idea underlying the above proof, namely to use a direct sum
  decomposition for obtaining information about invariants like the Hilbert
  series is very old and goes back at least to Riquier \cite{Riquier} and
  Janet \cite{Janet2} (in the context of partial differential equations).
  In fact, Janet gave already an explicit formula for the Hilbert function
  in terms of what is nowadays called a Janet basis.  Within algebra, it
  was mainly Stanley \cite{Stanley} who exploited much later the same idea.
  For this reason one often speaks of Stanley decompositions.  However, the
  special decompositions induced by Pommaret bases appeared already
  considerably earlier in the work of Rees \cite{Rees}.  These
  decompositions can also be used for the construction of resolutions, as
  Eliahou and Kervaire \cite{EK} showed first for the special case of
  stable ideals where they could obtain an explicit expression for the
  minimal resolution of the ideal.  As a by-product, they obtained this way
  via the Betti numbers \eqref{eq:hs} for this special case.  Later it was
  shown in \cite{wms:comb2} how their construction embeds into the theory
  of Pommaret bases and how it can consequently be generalised to
  quasi-stable ideals and (to some extent) to polynomial ideals in
  quasi-stable position.  However, in these more general situations one no
  longer obtains the minimal resolution and hence only upper bounds for the
  Betti numbers.  But as \eqref{eq:hs} is a simple consequence of the
  direct sum decomposition induced by the Pommaret basis independent of any
  explicit expression for the Betti numbers, it remains valid.
\end{remark}

Based on the proof above, we can derive an upper bound for the degree of
the numerator $N$ of the Hilbert series. Furthermore, we provide an
explicit formula for the coefficients of $N$ in terms of the Pommaret basis
$\H$.

\begin{proposition}
  \label{prop0sec8}
  Let $\HS_{\I}(t)=N(t)/(1-t)^D$ where
  $N(t)=a_0+a_1t+\cdots +a_\ell t^\ell$ with $a_{\ell}\ne 0$. If $\H$ is
  the Pommaret basis of $\I$, then $\deg(N)=\ell$ satisfies
  \begin{displaymath}
    \ell \le \max\{\deg(h)-|\mult{\X}{P}{h}|+D \ | \ h\in \H\}
  \end{displaymath}
  and the coefficients $a_{i}$ are given by
  \begin{equation}\label{eq:ai}
    \begin{aligned}
      a_i=&\binom{n-D+i-1}{n-D-1}-
           \sum_{\substack{h\in \H,\deg(h)\le i,\\ |\mult{\X}{P}{h}|\le D}}
           {(-1)^{i-\deg(h)}\binom{D-|\mult{\X}{P}{h}|}{i-\deg(h)}}\\[0.5\baselineskip]
          &{}-\sum_{\substack{h\in \H,\deg(h)\le i,\\ |\mult{\X}{P}{h}|> D}}
              {\binom{|\mult{\X}{P}{h}|-D+i-\deg(h)-1}
              {|\mult{\X}{P}{h}|-D-1}}\,.
    \end{aligned}
  \end{equation}
\end{proposition}

\begin{proof}
  By the proof of Thm.~\ref{th1sec8},
  $N(t)(1-t)^{n-D}=1-\sum_{j=1}^{s}{(1-t)^{n_j}t^{c_j}}$ and thus
  $\deg(N(t))\le \max\{n_j+c_j \ | \ j=1,\ldots ,s\}-n+D$ which immediately
  implies the bound for $\ell$.

  By the equality above,
  $N(t)=(1-t)^{D-n}-\sum_{j=1}^{s}{(1-t)^{n_j-n+D}t^{c_j}}$. Thus we have
  expressed $N$ as a sum of rational functions and now study the respective
  Taylor series.  The coefficient of $t^i$ in the series expansion of
  $(1-t)^{D-n}$ is equal to the number of terms of degree $i$ in $n-D$
  variables and hence to $\binom{n-D+i-1}{n-D-1}$, the first term in
  \eqref{eq:ai}. Now we must find the coefficient of $t^i$ in
  $(1-t)^{n_j-n+D}t^{c_j}$. Obviously, it vanishes for $c_j>i$. For
  $c_j\le i$ two cases arise. If $n_j-n+D\ge 0$ (or equivalently
  $|\mult{\X}{P}{h}|\le D$), then, by binomial expansion, the coefficient
  of $t^i$ is $(-1)^{i-c_j}\binom{n_j-n+D}{i-c_j}$ which yields the second
  summand in \eqref{eq:ai}. Otherwise, the coefficient of $t^i$ is equal to
  the number of terms of degree $i-c_j$ in $n-n_j-D$ variables and thus
  to $\binom{n-n_j-D+i-c_j-1}{n-n_j-D-1}$ leading to the last summand in
  \eqref{eq:ai}. \qed
\end{proof}

The above result leads to two simple corollaries relating the Hilbert
regularity with the parameters of the Pommaret basis and other invariants.

\begin{corollary}
  \label{prop1sec8}
  If $\H$ is the Pommaret basis of $\I$, then
  \begin{displaymath}
    \hilb(\I)\le
    \max\bigl\{0,\deg(h)-|\mult{\X}{P}{h}|+1 \mid h\in \H\bigr\}\,.
  \end{displaymath}
\end{corollary}

\begin{proof}
  Write $\HS_{\I}(t)=N(t)/(1-t)^D$ with $N(1)\ne 0$. It is well-known that
  $\hilb(\I)=\max\{0,\deg(N(t))-D+1\}$, see e.\,g.\
  \cite[Prop.~4.1.12]{Bruns}. Now Prop.~\ref{prop0sec8} immediately entails
  the assertion.  \qed
\end{proof}

\begin{corollary}\label{cor:hilb2}
  Assume that the Pommaret basis $\H$ of $\I$ contains a generator $h\in\H$
  having simultaneously the maximal degree and the maximal class among all
  elements of $\H$. Then
  \begin{itemize}
  \item[(1)] $\hilb(\I)=\max\{0,\deg(h)-|\mult{\X}{P}{h}|+1\}$,
  \item[(2)] $\deg(N(t))=\deg(h)-|\mult{\X}{P}{h}|+D$,
  \item[(3)] $\hilb(\I)+\depth(\I)=\max\{\depth(\I),\reg(\I)\}$.
  \end{itemize}
\end{corollary}

\begin{proof}
  We may assume that $h=h_s$. Since $n_s+c_s\ge n_i+c_i$ for all $i$, two
  first equalities follow from the proofs of Prop.~\ref{prop0sec8} and
  Cor.~\ref{prop1sec8}, respectively.  By \cite[Prop.~3.19]{wms:comb2} we
  have $\depth(\I)=|\mult{\X}{P}{h_s}|-1$ and $\reg(\I)=\deg(h_s)$ by
  \cite[Thm.~9.2]{wms:comb2} which entails the last assertion.  \qed
\end{proof}

\begin{remark}
  Mora \cite[Thm.~38.1.1]{Mora} claimed the equality
  $\reg(\I)=\hilb(\I)+\depth(\I)$ for arbitrary ideals.  However, it does
  not even necessarily hold for ideals in generic position\footnote{By
    generic position, we mean after a linear change of variables from a
    Zariski open set, see \cite{bayer_stillman} for more details.}. As a
  concrete counterexample, we consider the monomial ideal
  $\I=\langle x_1x_3, x_1x_2, x_1^2, x_2^3\rangle \subset \kk[x_1,x_2,x_3]$
  which is a generic initial ideal. Its generating set is already its
  Pommaret basis and therefore $\reg(\I)=3$.  The only term of maximal
  class is not of maximal degree and thus Cor.~\ref{cor:hilb2} cannot be
  invoked. , By \cite[Prop.~3.19]{wms:comb2}, $\depth(\I)=0$. On the other
  hand, $\HS_{\I}(t)=-(2t+1)/(-1+t)$ and therefore $\hilb(\I)=1$. Thus
  $\reg(\I)> \hilb(\I)+\depth(\I)$. Finally, we note that in general the
  inequality $\reg(\I)\ge \hilb(\I)+\depth(\I)$ does not hold. For example,
  consider the above ideal in the ring $\kk[x_1,\ldots ,x_7]$. Then,
  $\reg(\I)=3$, $\hilb(\I)=0$ and $\depth(\I)=4$.
\end{remark}

Finally, we provide an explicit expression for the degree of an ideal in
terms of its Pommaret basis.  Below, we denote by $f^{(i)}$ the $i$th
derivative of the function $f$.

\begin{theorem}
  \label{prop1.8}
  Let $\H$ be the Pommaret basis of $\I$. Then
  \begin{displaymath}
    \deg(\I)=\hspace{-0.5cm}
        \sum_{\substack{h\in \H,\\ 
              D \le |\mult{\X}{P}{h}| \le D+\deg(h)}}\hspace{-0.5cm}
        (-1)^{|\mult{\X}{P}{h}|-D+1}\binom{\deg(h)}{|\mult{\X}{P}{h}|-D}\,.
  \end{displaymath}
\end{theorem}

\begin{proof}
  We claim
  that
  \begin{displaymath}
    \deg(\I)=N(1)=
    \frac{(-1)^{n-D}}{(n-D)!}\Bigl(N(t)(1-t)^{n-D}\Bigr)^{(n-D)}|_{t=1}\,.
  \end{displaymath}
  Indeed, by the general Leibniz rule we have
  \begin{displaymath}
    \Bigl(N(t)(1-t)^{n-D}\Bigr)^{(n-D)}=
    \sum_{k_1+k_2=n-D}\frac{(n-D)!}{k_1!k_2!}N(t)^{(k_1)}\bigl((1-t)^{n-D}\bigr)^{(k_2)}\,.
  \end{displaymath}
  On the right hand side, all summands vanish at $t=1$ except for $k_1=0$
  and $k_2=n-D$ which proves the claim. 

  By the proof of Thm.~\ref{th1sec8},
  $N(t)(1-t)^{n-D}=1-\sum_{i=1}^{s}{(1-t)^{n_i}t^{c_i}}$. Thus there only
  remains to determine the derivatives of the right hand side. Applying
  again the general Leibniz rule, we obtain
  \begin{displaymath}
    \bigl((1-t)^{n_i}t^{c_i}\bigr)^{(n-D)}=
    \sum_{k_1+k_2=n-D}\frac{(n-D)!}{k_1!k_2!}
                    \bigl((1-t)^{n_i}\bigr)^{(k_1)}
                    (t^{c_i})^{(k_2)}\,.
  \end{displaymath}
  At $t=1$, all summands disappear except for $k_1=n_i\le n-D$ and
  $k_2\le c_i$. Such indices appear whenever $n_i+c_i\ge n-D$ and
  $n_i\le n-D$. By simple manipulations, one now obtains
  \begin{displaymath}
    \bigl((1-t)^{n_i}t^{c_i}\bigr)^{(n-D)}|_{t=1}=
    (-1)^{n_i}(n-D)!\binom{c_i}{n-D-n_i}
  \end{displaymath}
  which immediately yields the assertion. \qed
\end{proof}

\section{Dimension Depending Upper Bounds}
\label{sec:9}

We now exploit some of the results obtained in the last section to derive a
dimension depending upper bound for the degree of a homogeneous
ideal. Furthermore, as related subjects, we provide new dimension depending
upper bounds in the effective Nullstellensatz, in elimination theory and
for the ideal membership problem.  Let us quickly recall the used
notations.  $\I\subset \P$ is an ideal generated by the homogeneous
polynomials $f_1,\ldots ,f_k$.  If $\I$ is in quasi stable position, then
we denote its Pommaret basis by $\H$.  We write $d_i=\deg(f_i)$ and assume
that $d_1\ge \cdots \ge d_k$.  If an index $i>k$ appears, then we set
$d_{i}=1$.  Then the classical {\em B\'ezout bound} (see e.\,g.\
\cite[Thm.~1]{Heintz}, \cite[Lem.~3]{Brown2} or \cite[Lem.~2.95]{Ritscher})
asserts that $\deg(\I)\le d_1\cdots d_{\mu}$ with $\mu =\min\{k,n\}$.  We
will now improve this bound using the dimension\footnote{Although we are
  dealing with a homogeneous ideal, we will always work with the dimension
  as \emph{affine} variety.}  $D=\dim{\I}$.

\begin{lemma}\label{lem4}
  If the ideal $\I$ is in quasi stable position and $\H$ its Pommaret
  basis, then for each $i\le n$ the set $\H|_{x_{i}=\cdots=x_n=0}$ is the
  Pommaret basis of the ideal
  $\I|_{x_{i}=\cdots=x_n=0}\subseteq\kk[x_{1},\dots,x_{i-1}]$.
\end{lemma}

\begin{proof}
  Obviously, $\H|_{x_n=0}\subset \I|_{x_n=0}$ and no leading term of an
  element of $\H|_{x_n=0}$ is a Pommaret divisor of the leading term of
  another element.  It is well-known that the reverse lexicographic term
  order implies that, if $x_n$ divides the leading term of a polynomial,
  then it divides every term in the polynomial. It follows immediately that
  $\I|_{x_n=0}=\bigoplus_{h\in\H}\kk[\mult{\X}{P}{h|_{x_n=0}}]\cdot
  h|_{x_n=0}$ and $\H|_{x_n=0}$ is thus the Pommaret basis of
  $\I|_{x_n=0}$. For $i<n$ the claim is established using a simple
  induction. \qed
\end{proof}

\begin{proposition}\label{Cor}
  Suppose that the ideal $\I$ is in quasi stable position and of dimension
  $D>0$. Then $\deg(\I)=\deg(\I|_{x_{n-D+2}=\cdots=x_n=0})$.
\end{proposition}

\begin{proof}
  Let $\H$ be the Pommaret basis of $\I$.  It follows from
  Thm.~\ref{prop1.8} that generators with class greater than $n-D+1$ are
  not considered in the there provided formula for $\deg{(\I)}$.
  Together with Lem.~\ref{lem4}, this observation implies that the degrees
  of $\I$ and of $\I|_{x_{n-D+2}=\cdots=x_n=0}$ are identical. \qed
\end{proof}

\begin{example}
  Let
  $\I=\langle x_2x_3, x_1^2, x_1x_2x_4, x_2^3, x_1x_3^2x_4, x_1x_3^3,
  x_2^2x_4^2x_5, x_4^3x_2^2\rangle$ be an ideal in
  $\P=\kk[x_1,\ldots ,x_5]$. One can easily show that $\I$ is quasi stable
  and that $\dim(\I)=3$.  By Prop.~\ref{Cor},
  $\deg(\I)=\deg(\I|_{x_{4}=x_5=0})=1$, and this makes the computation of
  the degree of $\I$ less expensive. We note that generally Prop.~\ref{Cor}
  becomes false, if we set further variables to zero.  In our concrete
  example, one easily checks that the Hilbert series of
  $\I|_{x_3=x_{4}=x_5=0}$ is $t^3+2t^2+2t+1$ and therefore
  $\deg(\I|_{x_3=x_{4}=x_5=0})=6\neq\deg(\I)$.
\end{example}

\begin{corollary}
  \label{Bezout1}
  If $D>0$ then $\deg(\I)\le d_1\cdots d_{\mu}$ where
  $\mu=\min\{k,n-D+1\}$.
\end{corollary}

\begin{proof}
  By Rem.~\ref{linchen}, we may assume that $\I$ is in quasi stable
  position. By Prop.~\ref{Cor}, the degrees of $\I$ and of
  $\I|_{x_{n-D+2}=\cdots=x_n=0}$ are identical with the latter ideal lying
  in the ring $\kk[x_1,\ldots ,x_{n-D+1}]$.  Now the assertion follows by
  B\'ezout's bound. \qed
\end{proof}

\begin{theorem}[Dimension depending B\'ezout bound]
  \label{bezoutbound}
  If the ideal $\I$ has dimension $D$, then
  $\deg(\I)\le d_1\cdots d_{n-D}$.
\end{theorem}

\begin{proof}
  For $D=0$, this is just the classical B\'ezout bound. For $D>0$, we may
  assume that $\I$ is in quasi stable position by Rem.~\ref{linchen}.  By
  Prop.~\ref{Cor}, $\deg(\I)=\deg(\J)$ where
  $\J=\I|_{x_{n-D+2}=\cdots=x_n=0}\subseteq\kk[x_1,\ldots ,x_{n-D+1}]$ and
  it suffices to prove the desired upper bound for the latter ideal.  Since
  $\I$ is in quasi stable position and $D=\dim(\I)$,  a pure power of
  each variable $x_{1},\ldots ,x_{n-D}$ appears in $\LT(\I)$ and no pure
  power of $x_{n-D+1}$ belongs to $\LT(\I)$ (this follows e.\,g.\ from
  \cite[Prop.~3.15]{wms:comb2}). Therefore, $\dim(\J)=1$. 

  $\J$ is a homogeneous ideal generated by the polynomials
  $f_i|_{x_{n-D+2}=\cdots=x_n=0}$ with $i=1,\ldots ,k$. Since $\I$ is in
  quasi stable position, $\J$ is in quasi stable position, too,
  (Lem.~\ref{lem4}) and therefore $\J^{\sat}=\J:x_{n-D+1}^{\infty}$ (see
  e.g. \cite[Prop.~10.1]{wms:comb2}) . This
  implies that the degree of $\J^{\sat}$ equals the number of projective
  solutions (with multiplicity) of the system associated to $\J^{\sat}$
  (see e.\,g.\ \cite[Thm.~3.2]{Lazard81}).  Since the ideal $\J^{\sat}$ is
  saturated, this number equals the number of affine solutions (with
  multiplicity) of the system associated to the ideal $\J$ with
  $x_{n-D+1}=1$.  Obviously, the ideal $\J|_{x_{n-D+1}=1}$ is
  zero-dimensional and generated by the polynomials
  $f_i|_{x_{n-D+1}=1,x_{n-D+2}=\cdots=x_n=0}$ with $i=1,\ldots ,k$ (see
  \cite[Prop.~10.1]{wms:comb2}).  By B\'ezout's theorem, the number of
  solutions is thus bounded by $d_1\cdots d_{n-D}$ and hence
  $\deg(\J^{\sat})\le d_1\cdots d_{n-D}$. Our assertion now follows from
  the fact that $\deg(\J)=\deg(\J^{\sat})$.  \qed
\end{proof}

\begin{remark}
  An alternative proof of Thm.~\ref{bezoutbound} goes as follows.  Suppose
  that the ideal $\I$ is in quasi stable position. As we observed, for
  bounding $\deg(\I)$, we could only add the variables
  $x_{n-D+2},\ldots ,x_n$ into $\I$ and the main obstacle in estimating
  $\deg(\I)$ was the addition of $x_{n-D+1}$ into
  $\I|_{x_{n-D+2}=\cdots=x_n=0}$.  Using the notations in the proof above,
  we have $\deg(\I)=\deg(\J)=\deg(\J^{\sat})$ and
  $\J\subset \kk[x_1,\ldots ,x_{n-D+1}]$ is a one-dimensional ideal.  If
  $\H '=\H|_{x_{n-D+2}=\cdots=x_n=0}$ is the Pommaret basis of $\J$, then
  we set $\H' _1=\{h\in \H' \ | \ \cls(h)=n-D+1 \}$ and
  $\bar{\H' _1}=\{h/x_{n-D+1}^{\deg_{x_{n-D+1}}(h)} \ | \ h\in \H' _1\}$
  where $\deg_{x_{n-D+1}}(h)$ denotes the degree of $h$ in the variable
  $x_{n-D+1}$.  \cite[Prop.~10.1]{wms:comb2} asserts that
  $\H' \setminus \H' _1 \cup \bar{\H _1}$ is a (non reduced) Pommaret basis
  (and thus a Gr\"obner basis) of $\J^{\sat}$. It follows that $x_{n-D+1}$
  does not appear in $\LT(\J^{\sat})$ and therefore $\deg(\J)$ is equal to
  the degree of $\J^{\sat}|_{x_{n-D+1}=0} \subset \kk[x_1,\ldots ,x_{n-D}]$
  which is a zero-dimensional ideal.  On the other hand, we know that the
  ideal $\K$ generated by the polynomials
  $f_i|_{x_{n-D+1}=x_{n-D+2}=\cdots=x_n=0}$ with $i=1,\ldots ,k$ is
  zero-dimensional and $\K\subset \J^{\sat}|_{x_{n-D+1}=0}$. This yields
  that $\deg(\J^{\sat}|_{x_{n-D+1}=0})\le \deg(\K)$. So the desired
  inequality follows from B\'ezout's theorem.  For example, the generating
  set of the ideal
  $\I=\langle 5x_2^2+4x_2x_3+x_3^2, 2x_1x_2+x_1x_3+2x_2^2+x_2x_3,
  x_1^2+x_1x_2, x_1x_3^2+x_2x_3^2 \rangle $ in the polynomial ring
  $\kk[x_1,x_2 ,x_{3}]$ is already a Pommaret basis and $\dim(\I)=1$.  By
  the above argument, the set
  $\H'=\{5x_2^2, 2x_1x_2+2x_2^2, x_1^2+x_1x_2, x_1+x_2\}$ is a Gr\"obner
  basis for $\I^{\sat}$ and this observation entails $\deg(\I)=2$.
\end{remark}

\begin{remark}
  In general, this upper bound improves the upper bound $d_1^{n-D}$ for
  $\deg(\I)$ due to Masser and W\"ustholz \cite[Thm.~II]{Masser},
  cf.~\cite[Prop.~3.5]{bayer3}.  We also note that Lazard
  \cite[Prop.~1]{Daniel83} presented (without any proof or reference) a
  similar upper bound.
\end{remark}

\begin{example}
  Consider the homogeneous ideal
  \begin{multline*}
    \I=\langle x_5x_7+x_1x_8,\ x_6x_7+x_1x_9,\ x_6x_8+x_5x_9,\ x_5x_2+x_1x_3,\\
    x_6x_2+x_1x_4,\ x_6x_3+x_5x_4,\ x_8x_2+x_7x_3,\ x_9x_2+x_7x_4,\ x_9x_3+x_8x_4
    \rangle
  \end{multline*}
  in the polynomial ring $\kk[x_1, x_2, x_3, x_4, x_5, x_6, x_7, x_8, x_9]$
  appearing in the work of Eisenbud and Sturmfels \cite{Eisenbud}.  We note
  that $\dim(\I)=3$, $\deg(\I)=24$ and $\I$ is generated by $9$ quadratic
  polynomials.  The classical B\'ezout bound yields the estimate
  $\deg(\I)\le 2^9=512$, while Thm.~\ref{bezoutbound} says
  $\deg(\I)\le 2^6=64$.
\end{example}

\begin{example}
  The so-called Mora-Lazard-Masser-Philippon-Koll\'ar example \cite{Brown}
  shows that the degree bound of Thm.~\ref{bezoutbound} is sharp.  For any
  sequence of degrees $d_{1},\ldots ,d_{n-1}$, let $\I$ by the ideal
  generated by the set
  \begin{displaymath}
    A=\bigl\{x_{1}^{d_{1}}-x_2x_{n}^{d_1-1},\ x_{2}^{d_{2}}-x_3x_{n}^{d_2-1},\ldots
    ,\ x_{n-2}^{d_{n-2}}-x_{n-1}x_{n}^{d_{n-2}-1},\ x_{n-1}^{d_{n-1}}\bigr\}\,.
  \end{displaymath}
  The first Buchberger criterion (see e.\,g.\ \cite{little}) shows easily
  that $A$ is the reduced Gr\"obner basis of $\I$. Therefore,
  $\LT(\I)=\langle x_{1}^{d_{1}},\ldots ,x_{n-2}^{d_{n-2}},
  x_{n-1}^{d_{n-1}} \rangle$ entailing that
  $\HS_{\I}(t)=(1-t^{d_1})\cdots (1-t^{d_{n-1}})/(1-t)^{n}$ and thus
  $\deg(\I)=d_1\cdots d_{n-1}$.
\end{example}

We finally discuss some new dimension depending upper bounds for the
effective Nullstellensatz, for elimination theory and for the ideal
membership problem.  We first briefly review some known results related to
the effective Nullstellensatz that we will use in the rest of this
section. For a sequence $d_1\ge \cdots \ge d_k$ of positive integers, let
\begin{equation*}
 \N(d_1,\ldots ,d_k;n)=
 \begin{cases}
   {\prod_{i=1}^{k}d_i} & \text{if $n\ge k\ge 1$} \\
   {d_k\prod_{i=1}^{n-1}d_i} & \text{if $k>n>1$} \\
   {d_k} & \text{if $n=1.$} \\
 \end{cases}
\end{equation*}
In an effective Nullstellensatz, one considers statements as follows.  If a
homogeneous polynomial $f$ belongs to the radical of the ideal generated by
the polynomials $f_i$, then there exists a positive integer $e$ and
polynomials $g_1,\ldots ,g_k\in \P$ such that $f^e=g_1f_1+\cdots +g_kf_k$
with $e\le \N(d_1,\ldots ,d_k;n)$ and
$\deg(g_if_i)\le \deg(f)\N(d_1,\ldots ,d_k;n)$.  This result also holds for
non-homogeneous polynomials, if we replace the factor $\deg(f)$ by
$1+\deg(f)$.  The smallest integer $e$ such that for every polynomial
$f\in \sqrt{\I}$ we have $f^e\in \I$ is called the \emph{N\oe ther
  exponent} of $\I$ and it is denoted by ${\rm e}(\I)$.

Let $X$ be a variety of dimension $n$ and $\deg(X)$ its degree. Let
$f_1,\ldots ,f_k\in \kk[X]$ be a sequence of polynomials so that
$d_i=\deg(f_i)$ with $d_1\ge \cdots \ge d_k$. Then, the above assertion
holds true for the $f_i$'s if in the above upper bound we replace
$\N(d_1,\ldots ,d_k;n)$ by $\deg(X)\N(d_1,\ldots ,d_k;n)$. Thus,
Thm.~\ref{bezoutbound} could be useful for finding an upper bound for
$\deg(X)$, if it is unknown. For further details on this topic, we refer to
e.g.\ \cite{Kollar,Sombra,Jelonek}.

We are now concerned with obtaining an upper bound for the N\oe ther
exponent $\mathrm{e}(\I)$. For unmixed\footnote{An ideal is called unmixed
  if all its associated prime ideals have the same dimension.} ideals, it has
been shown that ${\rm e}(\I)\le d_1\cdots d_{n-D}$
cf.~\cite[Lem.~4]{Brown2}.  However, this inequality does not hold in
general.

\begin{example}\label{exconj}
  Consider the ideal\footnote{This example has been provided by David
    Masser (private communication).}
  $\I=\langle x^4+x^3y, x^3+y^3 \rangle \subset \kk[x,y]$.  We have
  $\I=\langle x+y \rangle \cap \langle x^3,y^3\rangle, \dim(\I)=1$ and
  $d_1=4$. On the other hand, we note that $x+y\in \sqrt{I}$ however
  $(x+y)^4\notin \I$ and $(x+y)^5\in \I$. Therefore, $e\geq5\not\le d_1=4$.
\end{example}

In characteristic zero, we now provide a refinement of this bound for all
one-dimensional ideals using the next theorem due to Lazard
\cite[Thm.~2]{Daniel83}.

\begin{theorem}\label{laz3}
  Let the ideal $\I$ be in generic position and assume that $D\le 1$. Then
  $\deg{(\I,\prec)}\le d_1+\cdots +d_r-r+1$ where $r=n-\depth(\I)$.
\end{theorem}

\begin{proposition}\label{prop:edim1}
  Let $\I$ be an ideal of dimension $D=1$ over a field of characteristic
  zero. Then
  ${\rm e}(\I) \le \max{\{ d_1\cdots d_{n-1}, d_1+\cdots +d_r-r+1\}}$ where
  $r=n-\depth(\I)$.
\end{proposition}

\begin{proof}
  W.l.o.g.\ we may assume that $\I$ is in generic position.  Since $D=1$
  and $\I$ is homogeneous, a primary decomposition
  $\I=Q_1\cap \cdots \cap Q_s\cap Q$ exists where each $Q_i$ is a
  $\mathfrak{p}_i$-primary component of $\I$ of dimension $1$ and $Q$ is
  $\mf$-primary.  Let $\ell_i$ be the length of $Q_i$. By
  Rem.~\ref{rem:lendeg},
  $\deg(\I)=\sum_{i=1}^{s}{\ell_i\deg(\mathfrak{p}_i)}$. We note that
  $\sqrt{I}=\mathfrak{p}_1\cap \cdots \cap \mathfrak{p}_s$.  For any
  polynomial $f\in \mathfrak{p}_1\cap \cdots \cap \mathfrak{p}_s$, one has
  $f^{\deg(\I)}\in Q_1\cap \cdots \cap Q_s$ (cf.\
  Rem.~\ref{rem:lendeg}). On the other hand, one easily sees that any
  $f\in \mf$ satisfies $f^{\sat(\I)}\in Q$. Hence we find
  ${\rm e}(\I)\le \max \{\deg(\I),\sat(\I)\}$. Now the assertion follows
  from Thm.~\ref{laz3} and the fact that in generic (more precisely, in
  stable) position $\sat(\I)\le \reg(\I)=\deg(\I,\prec)$ (see e.\,g.\
  \cite[Thms.~5.5.7 and 5.5.15, Prop.~5.5.28]{wms:invol}). \qed
\end{proof}

\begin{example}
  We consider again the ideal of Example \ref{exconj}. There $\depth(\I)=0$
  and by Prop.~\ref{prop:edim1} we get
  ${\rm e}(\I)=5\le \max{\{4,4+3-2+1\}}=6$.
\end{example}

\begin{remark}
  The restriction to characteristic zero is required in the last step of
  the proof of Prop.~\ref{prop:edim1}.  Only in characteristic zero we can
  always achieve stable position by a linear transformation.  In positive
  characteristic $p>0$, the ideal
  $\I=\langle x_{1}^{p}, x_{2}^{p}\rangle\subset\kk[x_{1},x_{2}]$ is not in
  stable position and invariant under any linear transformation. For such
  ideals we only have the inequality $\reg{(\I)}\geq\deg{(\I,\prec)}$ and
  cannot conclude any relationship between $\sat{(\I)}$ and
  $\deg{(\I,\prec)}$.  For the above ideal
  $\sat{(\I)}=\reg{(\I)}=2p-1>\deg{(\I,\prec)}=p$.  We get the same values
  for these invariants, if we consider $\I$ as an ideal in
  $\kk[x_{1},x_{2},x_{3}]$ so that it is one-dimensional.
  Prop.~\ref{prop:edim1} remains correct for this example, as the bound
  comes now from the first term which is $p^{2}$.
\end{remark}

In the remainder of this section, unless explicitly stated otherwise, we
skip the assumption that we are dealing with homogeneous polynomials and
ideals, but keep otherwise our notations.  We present first a dimension
depending bound for the representation problem related to N\oe ther
normalization \cite{Alicia} (see also \cite{Kartzer,MR}).  Recall that an
ideal $\I\subset \P$ is in \emph{N\oe ther position}, if the ring extension
$K[x_{n-D+1} ,\ldots ,x_{n} ] \hookrightarrow \P/\I$ is integral
\cite{Eisenbud_book}.  In this case, for each index $1\leq i\leq n-D$, the
intersection $\I\cap \kk[x_i,x_{n-D+1},\ldots ,x_n]$ is non-empty and, by
\cite[Prop.~1.7]{Alicia}, contains a \emph{witness polynomial} $h_i$ which
is monic in $x_i$ and can be represented in the form
$h_i=g_1f_1+\cdots +g_kf_k$ with coefficients $g_1,\ldots ,g_k\in \P$ such
that $\deg(g_jf_j)\le d^n(d^n+1)$ for $d=d_1$. In this estimate, the first
factor $d^n$ represents an upper bound for $\N(d_1,\ldots ,d_k;n)$.  Mayr
and Ritscher \cite{MR} proved the following improvement.

\begin{proposition}(\cite[Thm.~10]{MR})\label{MR} 
  In the above described situation, the coefficients $g_{j}$ can be chosen
  such $\deg(g_jf_j)\le (d_1\cdots d_{n-D})^2$.
\end{proposition}

We will now improve this bound to $3d_1\cdots d_{n-D}$ and study its
application in the membership problem.  We denote by $\leftexp{h}{\P}$ the
ring $\kk[x_1,\ldots ,x_{n+1}]$ where $x_{n+1}$ is a new variable. For any
polynomial $f\in \P$, we consider its homogenization
$\leftexp{h}{f}=x_{n+1}^{\deg(f)}f(x_1/x_{n+1},\ldots ,x_n/x_{n+1})\in
\leftexp{h}{\P}$.  For an ideal $\I\subset \P$, its homogenization is
defined as
$\leftexp{h}{\I}=\langle\leftexp{h}{f} \ | \ f\in \I\rangle \subset
\leftexp{h}{\P}$. We need the following result due to Sombra
\cite[Lem.~3.15]{Sombra} which allows us to extract a regular sequence
which remains regular after homogenization. We denote by $\I_i$, for
$i=1,\ldots ,k$, the ideal generated by $f_1,\ldots ,f_i$ with the
convention $\I_0=\li 0\ri$.

\begin{lemma}(\cite[Lem.~3.15]{Sombra})\label{regular} Let
  $f_1,\ldots ,f_k$ be a regular sequence. Then polynomials
  $p_1,\ldots ,p_k$ and $q_1,\ldots , q_k$ exist in $\P$ such that for any
  $i$
  \begin{itemize}
  \item[$\bullet$] $\leftexp{h}{p_i}=x_{n+1}^{c_i}\leftexp{h}{f_i}+q_i$
    with $q_i\in \leftexp{h}{\I_{i-1}}$ and
    $c_i\le \max \{\deg(\leftexp{h}{\I_{i-1}}),\deg(f_i)\}$
  \item[$\bullet$]
    $\deg(p_i)\le \max \{\deg(\leftexp{h}{\I_{i-1}}),\deg(f_i)\}$
  \item[$\bullet$] $\leftexp{h}{p_1},\ldots ,\leftexp{h}{p_k}$ forms a
    regular sequence in $\leftexp{h}{\P}$.
  \end{itemize}
\end{lemma}

We remark that Sombra \cite{Sombra} assumed the conditions
$d_2\ge \cdots \ge d_k\ge d_1$ which may make our next bounds sharper,
however, using the fact that any permutation of a regular sequence is a
regular sequence and for simplicity we continue with our restrictions on
the degrees. We also use the following two results related to regular
sequences.  The next proposition can be found e.\,g.\ in
\cite[Lem.~2.81]{Ritscher} or \cite[Lem.~9]{MR}.  However, we prove it for
the sake of completeness.  In the proof, we apply the well-known fact that
the set of all zero-divisors for an ideal in $\P$ is the union of all its
minimal prime ideals.

\begin{proposition}\label{Sec3:Lem1}
  If the field $\kk$ is infinite, then there are a strictly decreasing
  sequence of integers $1\le j_{n-D}<\cdots <j_{1}=k$ and homogeneous
  polynomials $h_{i,j}\in \P$ such that each of the polynomials
  $g_i=f_i+h_{i,i+1}f_{i+1}+\cdots+h_{i,k}f_{k}$ for
  $i=j_1,\ldots ,j_{n-D}$ is homogeneous of degree $d_i$ and such that
  $g_{j_1},\ldots,g_{j_{n-D}}$ form a regular sequence in $\P$.
\end{proposition}

\begin{proof}
  We follow the proof given by Ritscher \cite[Lem.~2.81]{Ritscher} filling
  in some missing details.  We show by induction that for each
  $1\le r\le n-D$ there exists a regular sequence $g_{j_1},\ldots ,g_{j_r}$
  and $\dim(\langle f_{j_r},\ldots ,f_k\rangle)=n-r$.  For the base case,
  we take $j_1=k$ and $g_{j_1}=f_k\ne 0$.  Obviously,
  $\dim{(\langle f_k\rangle)}=n-1$.  For the inductive step, assume that
  $g_{j_1},\ldots,g_{j_{r}}$ for $r<n-D$ is a regular sequence in $\P$ and
  $\dim{(\langle f_{j_r},\ldots ,f_k\rangle)}=n-r$.  Let
  $P=\{\mathfrak{p}_1,\ldots ,\mathfrak{p}_t\}$ be the set of all
  associated primes of $\J=\langle g_{j_1},\ldots,g_{j_{r}} \rangle$.
  Then, by Macaulay's Unmixedness theorem (see e.\,g.\
  \cite[Thm.~2.1.6]{Bruns}), we have $\dim(\mathfrak{p}_i)=n-r$ for each
  $i$.  Since $\dim(\I)=n-D$ and $n-r>D$, there exists an integer $j_{r+1}$
  such that
  $\dim{(\langle f_{j_{r+1}},\ldots ,f_k\rangle)}\le n-(r+1)\ge D$.  Let
  $j_{r+1}$ be the maximum integer with this property.

  Consider the $\kk$-linear space
  $S=\prod_{i=j_{r+1}}^{k}{\kk^{\binom{d_{j_{r+1}}-d_i+n-1}{n-1}}}$ and for
  each integer $\ell=1,\ldots ,t$ the subspace
  \begin{displaymath}
    S_{\ell}=\left\{(a_{i,\alpha})_{j_{r+1}\le i \le k,
        |\alpha|=d_{j_{r+1}}-d_i}\in S \ | \
      \sum_{i=j_{r+1}}^{k}{\sum_{|\alpha|=d_{j_{r+1}}-d_i}
        {a_{i,\alpha}x^{\alpha}f_{i}}}\in\mathfrak{p}_\ell \right\}.
  \end{displaymath}
  We claim that $S_\ell$ is a {\em proper} subspace of $S$ for each $\ell$.
  For a proof by reductio ad absurdum assume that $S_\ell=S$ for some
  $\ell$.  Then we have $x_j^{d_{j_{r+1}}-d_i}f_{i}\in \mathfrak{p}_\ell$
  for each $j$ and for each $i=j_r+1,\ldots ,j_{r+1}$ and hence
  $f_{j_{r+1}},\ldots ,f_{k}\in \mathfrak{p}_\ell$ since
  $\mathfrak{p}_\ell$ is a prime ideal.  By construction, the ideal
  generated by $f_{j_{r+1}},\ldots ,f_{k}$ is of dimension $n-(r+1)$ which
  yields a contradiction.  

  Since $\kk$ is assumed to be infinite, $S\ne S_1\cup \cdots \cup S_t$ by
  elementary linear algebra.  Choose a tuple
  $(a_{i,\alpha})_{i,\alpha}\in S\setminus (S_1\cup \cdots \cup S_t)$.
  Then the corresponding polynomial
  $g'_{j_{r+1}}=\sum_{i=j_{r+1}}^{k}
  {\sum_{|\alpha|=d_{j_{r+1}}-d_i}{a_{i,\alpha}x^{\alpha}f_{i}}}$ is a
  non-zero divisor on $\P/\J$.  Note that here $f_{j_{r+1}}$ is multiplied
  only by a constant.  We show now that this constant does not vanish.
  Indeed, otherwise a linear combination of the polynomials
  $f_{j_{r+1}+1},\ldots ,f_{k}$ was a non-zero divisor on $\P/\J$ implying
  that the depth of the ideal $\langle f_{j_{r+1}+1},\ldots ,f_{k} \rangle$
  was greater than its dimension $n-r$ which is not possible.  Finally,
  dividing $g'_{j_{r+1}}$ by the coefficient of $f_{j_{r+1}}$ yields a new
  polynomial $g_{j_{r+1}}$ of the desired form to extend our regular
  sequence. \qed
\end{proof}

We note that one obtains a regular sequence for a generic choice of the
polynomials $h_{i,j}$.  Furthermore, this proposition implies that for a
given ideal $\I=\langle f_1,\ldots, f_k \rangle$ we may assume w.l.o.g.\
that $f_1,\ldots,f_{n-D}$ is a regular sequence.  Since we have
$\deg(g_i)\le \deg(f_i)$ for each $i$, this assumption may only increase
the following upper bounds.

\begin{proposition}(\cite[Cor.~3.5, page 107]{Monique})\label{Monique107}
  Let $f_1,\ldots, f_k$ be a regular sequence of homogeneous polynomials
  and assume that $\I=\langle f_{1},\dots,f_{k}\rangle$ is in N\oe ther
  position. Then, $\deg(\I,\prec)\le d_1+\cdots +d_k-k+1$.
\end{proposition}

\begin{theorem}\label{thmnoether}
  If $d_k\ge 2$, then a linear change of variables $\phi$ exists such that
  the transformed ideal $\phi(\I)$ is in N\oe ther position.  Furthermore,
  there are polynomials $h_i\in \I\cap \kk[x_i,x_{n-D+1},\ldots ,x_n]$ for
  $i=1,\ldots ,n-D$ and coefficients $g_1,\ldots ,g_k\in \P$ such that
  $h_i$ is monic in $x_i$ and $h_i=g_1f_1+\cdots +g_kf_k$ with
  $\deg(g_jf_j)\le 3d_1\cdots d_{n-D}$.
\end{theorem}

\begin{proof}
  By Prop.~\ref{Sec3:Lem1}, we may assume that the sequence
  $f_1,\ldots ,f_{n-D}$ is regular. Applying Lem.~\ref{regular} provides
  then polynomials $p_1,\ldots ,p_{n-D}$ such that their homogenizations
  form a regular sequence in $\leftexp{h}{\P}$ and satisfy
  $\deg(\leftexp{h}{p_j})\le \max
  \{\deg(\leftexp{h}{\I_{j-1}}),\deg(f_j)\}$.  $\leftexp{h}{\I_{j-1}}$ is
  the saturation of the ideal generated by
  $\leftexp{h}{f_1},\ldots, \leftexp{h}{f_{j-1}}$ w.r.t. $x_{n+1}$ and thus
  $\deg(\leftexp{h}{\I_{j-1}})\le \deg(\langle \leftexp{h}{f_1},\ldots,
  \leftexp{h}{f_{j-1}} \rangle )\le d_1\cdots d_{j-1}$. This observation
  implies that $\deg(\leftexp{h}{p_1})= d_1$ and
  $\deg(\leftexp{h}{p_j})\le d_1\cdots d_{j-1}$ for each $j>1$. Let
  $\J\subset \leftexp{h}{\P}$ be the homogeneous ideal generated by the
  polynomials $\leftexp{h}{p_j}$.  The N\oe ther normalization lemma
  asserts the existence of a linear change of variables $\phi$ such that
  $\phi(\J)$ is in N\oe ther position.  Since at the end, we will set
  $x_{n+1}=1$, we choose $\phi$ by ignoring $x_{n+1}$ in the N\oe ther
  normalization process which is always possible. It is easy to see that
  the sequence $\phi(\leftexp{h}{p_1}),\ldots , \phi(\leftexp{h}{p_{n-D}})$
  remains regular.

  We consider now a degree reverse lexicographic order with
  $x_{n+1}\prec x_{n}\prec \cdots \prec x_{n-D+1}\prec x_i$ and
  $x_i\prec x_j$ for $j\ne i$ and $j<n-D$.  Since $\phi(\J)$ is in N\oe
  ther position, the reduced Gr\"obner basis $G$ of $\phi(\J)$ for $\prec$
  contains a polynomial $w_i\in \kk[x_i,x_{n-D+1},\ldots ,x_n,x_{n+1}]$
  which is monic in $x_i$. Since $\phi(\J)$ is in addition generated by a
  regular sequence, Prop.~\ref{Monique107} implies that the degrees of the
  elements of $G$ is at most
  $\sum_{j=1}^{n-D}{(\deg(\leftexp{h}{p_j})-1)+1}\le d_1+d_1+d_1d_2+\cdots
  +d_1\cdots d_{n-D-1}$. Using a simple induction and the fact that
  $2\le d_j$ for all $j$, we conclude that
  $d_1+d_1+d_1d_2+\cdots +d_1\cdots d_{n-D-1}\le d_1\cdots d_{n-D}$.  Now,
  there exist coefficients $a_j\in \leftexp{h}{\P}$ such that
  $w_i=a_1\leftexp{h}{p_1}+\cdots +a_{n-D}\leftexp{h}{p_{n-D}}$.  We note
  that $ w_i|_{x_{n+1}=1}\in \phi(\I)$ is a monic polynomial in $x_i$ with
  coefficients in $\kk[x_{n-D+1},\ldots ,x_n]$.  However, there still
  remains to find the maximum degree of the representation of this
  polynomial in terms of the original generators $f_i$.

  Thus, we next aim at expressing the polynomials $\leftexp{h}{p_j}$ in
  terms of the generators $\leftexp{h}{f_1},\ldots,
  \leftexp{h}{f_{n-D}}$. We know that
  $\leftexp{h}{p_j}=x_{n+1}^{c_j}\leftexp{h}{f_j}+q_j$ with
  $q_j\in \leftexp{h}{\I_{j-1}}$ and $\deg(q_j)\le d_1\cdots d_{j-1}$. From
  \cite[Lem.~3.18]{Sombra}, we deduce that for some exponent
  $\mu_j\le 2d_1\cdots d_{j-1}$ the product $x_{n+1}^{\mu_j}q_j$ belongs to
  the ideal generated by $\leftexp{h}{f_1},\ldots,
  \leftexp{h}{f_{j-1}}$. This shows that $x_{n+1}^\mu w_i$ with
  $\mu =2d_1\cdots d_{n-D}$ can be written as a linear combination of the
  polynomials $\leftexp{h}{f_j}$ and the maximal degree of this expression
  is at most $3d_1\cdots d_{n-D}$.  If we set
  $h_i=x_{n+1}^\mu w_i|_{x_{n+1}=1}$, then the desired conditions hold for
  $h_i$, and this terminates the proof. \qed
\end{proof}

\begin{remark}
  The new bound given by the above theorem improves the existing results
  including the bound $d^n(d^{n}+1)$ stated in \cite[Sec.~1]{Alicia} where
  $d=\max\{d_1,\ldots ,d_k\}$ and also the bound $(d_1\cdots d_{n-D})^2$
  from \cite[Thm.~10]{MR}.
\end{remark}

We now state some consequences of this theorem.  Dickenstein et
al.~\cite{Alicia} applied their variant of Thm.~\ref{thmnoether} to give a
bound for the degree w.r.t.\ only a subset of variables for the membership
problem. Following their proof and using our new bound, we obtain the next
result.

\begin{proposition}\label{propelim}
  Assume $d_k\ge 2$.  A polynomial $f\in \P$ lies in the ideal $\I$, if and
  only if coefficients $g_1,\ldots ,g_k\in \P$ exist such that
  $f=g_1f_1+\cdots +g_kf_k$ and the degree of each summand $g_if_i$ w.r.t.\
  the variables $x_1,\ldots ,x_{n-D}$ is at most
  $\max\{\deg(f), d_1+3(n-D)d_1\cdots d_{n-D}\}+3d_1\cdots d_{n-D}$.
\end{proposition}

Consider the element $f=g_1f_1+\cdots +g_kf_k$ in the ideal generated by
the polynomials $f_1,\ldots ,f_k$ of degrees $d_1\ge \cdots \ge d_k$. The
first upper bound for the coefficients,
$\deg(g_i)\le \deg(f)+2(kd_1)^{2^{n-1}}$, was established by Hermann
\cite{Her}.  We will mimic the proof of \cite[Thm.~5]{Kartzer} to give a
sharper upper bound.  We first recall a generic degree bound due to
Hermann, see e.g.~\cite[page~312]{MM}.

\begin{proposition}\label{MM}
  Consider a linear system of equations
  \begin{displaymath}
    \sum_{j=1}^{s}{h_{ij}}X_j=c_i, \ \ \ \ i=1,\ldots ,t
  \end{displaymath}
  with coefficients $h_{ij},c_i\in \P$ which has at least one solution.
  Then, the system possesses a solution $g_1,\ldots ,g_s$ with
  $\deg(g_i)\le c+(ds)^{2^n}$ where $d=\max\{\deg(h_{ij})\}$ and
  $c=\max\{\deg(c_i)\}$.
\end{proposition}

\begin{theorem}
  Assume $d_k\ge 2$.  A polynomial $f\in \P$ lies in $\I$, if and only if
  coefficients $g_1,\ldots ,g_k\in \P$ exits such that
  $f=g_1f_1+\cdots +g_kf_k$ and the degree of each $g_i$ is at most
  $\deg(f)+(kd_1^{D})^{2^{n-D}}$.
\end{theorem}

\begin{proof}
  W.l.o.g., we may assume that $\I$ is in N\oe ther position.  Any ideal
  member $f\in \I$ can be written as a linear combination
  $f=g_1f_1+\cdots +g_kf_k$ where the degree of each summand $g_if_i$
  w.r.t.  $x_1,\ldots ,x_{n-D}$ is at most $B$ with $B$ the bound in
  Prop.~\ref{propelim}.  Set $\P_D=\kk[x_{n-D+1},\ldots ,x_n]$. We consider
  now $g_i,f_i$ and $f$ as elements of the polynomial ring
  $\P_D[x_1,\ldots ,x_{n-D}]$. This leads to representations
  $g_i=\sum_{j}{u_{ij}m_{ij}}$, $f_i=\sum_{j}{w_{ij}m'_{ij}}$ and
  $f=\sum_{j}{v_{j}m''_{j}}$ with coefficients $u_{ij},w_{ij},v_j\in \P_D$
  and terms $m_{ij},m'_{ij},m''_j\in \kk[x_1,\ldots ,x_{n-D}]$.  These
  satisfy for each $i,j$ the estimates $\deg(m_{ij}m'_{ij})\le B$,
  $\deg(w_{ij})\le d_1$ and $\deg(v_{j})\le \deg(f)$.  

  Since we look for an upper bound for the degrees of the $g_i$, we
  consider the coefficients $u_{ij}$ as unknowns over $\P_D$ and try to
  bound their degrees using Prop.~\ref{MM}.  We enter the above
  representations of $g_i,f_i$ and $f$ into the linear combination
  $f=g_1f_1+\cdots +g_kf_k$ and extract linear equations for the $u_{ij}$
  over $\P_D$. By equating the coefficients of each term (in the
  variables $x_1,\ldots ,x_{n-D}$) of degree at most $B$ on both sides, we
  derive a linear system of equations over $\P_D$. Solving the resulting
  system yields the coefficients $g_1,\ldots ,g_k$. The number of variables
  in each linear equation is at most $k$ times the number of terms in
  $x_{n-D+1},\ldots ,x_{n}$ of degree at most $d_1$, i.\,e.\ at most
  $k\binom{D+d_1}{D}\le kd_1^{D}$.  Moreover, the right hand side in each
  linear equation is a polynomial of degree at most
  $\deg(f)$. Prop.~\ref{MM} implies now the existence of a solution with
  $\deg(g_i)\le \deg(f)+(kd_1^{D})^{2^{n-D}}$ and this yields the desired
  bound. \qed
\end{proof}

We conclude this section by providing another consequence of
Thm.~\ref{thmnoether}.  Mayr and Ritscher \cite[Thm.~36]{MR} applied
Prop.~\ref{MR} to prove the upper bound
$$2\left(\frac{1}{2}((d_1\cdots d_{n-D})^{2(n-D)}+d_1) \right)^{2^D}$$
for the degrees of the elements of any reduced Gr\"obner basis of $\I$. In
their proof, they exploited that the homogenization of $\I$ contains a
homogeneous regular sequence of degree at most $(d_1\cdots d_{n-D})^2$. We
now improve this result.  The proof of Thm.~\ref{thmnoether} entails that
there are polynomials $p_1,\ldots ,p_{n-D}\in \I$ such that
$\leftexp{h}{p_1},\ldots ,\leftexp{h}{p_{n-D}}$ is a regular sequence of
degree at most $d_1\cdots d_{n-D}$, cf.~\cite[Lem.~35]{MR}, which yields
the following sharper bound.

\begin{corollary}
  If $d_k\ge 2$, then the degrees of the elements of any reduced Gr\"obner
  basis of $\I$ are bounded by
  \begin{displaymath}
    2\left(\frac{1}{2}((d_1\cdots d_{n-D})^{n-D}+d_1) \right)^{2^D}\,.
  \end{displaymath}
\end{corollary}

\section*{Acknowledgments.} 

The authors would like to thank the anonymous reviewers for their helpful comments which helped us to improve the manuscript. The third author received funding from the European Union's Horizon 2020 research and innovation programme under grant agreement No H2020-FETOPEN-2015-CSA 712689.




\end{document}